\documentclass[
final
]{dmtcs-episciences}

\usepackage[utf8]{inputenc}
\usepackage{subfigure}

\usepackage[round]{natbib}

\newtheorem{thm1}{Theorem}

\newtheorem{cor1}{Corollary}

\newtheorem{prop1}{Proposition}

\author{Danilo Kor\v ze \affiliationmark{1}\thanks{Supported  by the Slovenian Research Agency under the grant J2-7357.}
  \and Aleksander Vesel \affiliationmark{2}\thanks{Supported
  by the Slovenian Research Agency under the grants P1-0297, J1-7110 and  J1-9109.}}
\title{Packing coloring of  generalized Sierpi\'nski graphs}
\affiliation{
Faculty of Electrical Engineering and Computer Science, University of Maribor,  Maribor, Slovenia \\
Faculty of Natural Sciences and Mathematics, University of Maribor, Maribor, Slovenia}

\received{2018-10-1}
\revised{2019-1-22}
\accepted{2019-1-24}

\begin{document}
\publicationdetails{21}{2019}{3}{7}{4862}
\maketitle
\begin{abstract}
  The packing chromatic number $\chi_{\rho}(G)$ of a graph $G$ is the smallest
integer $c$ such that the vertex set $V(G)$ can be partitioned into 
sets $X_1, . . . , X_c$, with the condition that vertices in $X_i$ have pairwise distance greater than $i$.
In this paper, we consider the packing chromatic number of several families of Sierpi\'nski--type graphs.
We establish  the packing chromatic numbers of  
generalized  Sierpi\'nski graphs $S^n_G$ where $G$ is a path or a cycle (with exception of the cycle of length five) 
as well as for two families where $G$ is a connected graph of order four. 
Furthermore, we prove that the packing chromatic number in the family
of Sierpi\'nski--triangle graphs $ST_4^n$  is bounded from above by 20. 
\end{abstract}

\section{Introduction}


A  {\em $c$-coloring} of a graph $G$ is a function $f$ from  $V(G)$ onto a set $ C = \{1, 2, \ldots, c \}$ (with no additional constraints). The elements of $C$ are called {\em colors}, while the set of vertices with the image (color) $i$ 
is denoted by $X_i$.   
Let $u,v$ be vertices of a graph $G$. The {\em distance} between $u$ and $v$ in $G$, denoted by $d_G(u, v)$,
  equals the length of a shortest $u,v$-path (i.e. a path between $u$ and $v$)  
in $G$. 

Let $f$ be a $c$-coloring of a graph $G$ with the corresponding sequence of color classes  $X_1, . . . , X_c$.
 If each color class $X_i$ is a set of vertices with the property that 
any distinct pair $u, v \in  X_i$ satisfies $d_G(u, v) > i$,
then $X_i$ is said to be an {\em $i$-packing}, while the sequence $X_1, . . . , X_c$ is called a {\em packing $c$-coloring}.
The smallest integer $c$ for which there exists a packing $c$-coloring of $G$
is called the {\em packing chromatic number of $G$} and it is denoted by  $\chi_{\rho}(G)$, see \cite{BrKlRa,Goddard}.

If $n$ is a natural number, let $[n]$ denote the set $\{0,1,\ldots,n-1\}$,

Let $G$ be an undirected graph with vertex set  $[k]$. The {\em generalized Sierpi\'nski
graph $S^n_G$ of $G$ of dimension $n$} is the graph with vertex set $[k]^n$,  while 
vertices $u, v \in V (S^n_G)$  are adjacent if and only if  there exists $i \in \{1, 2, 3, \ldots, n\}$ such that:

(i) $u_j = v_j$ if $j < i$,

(ii) $u_i  \not =  v_i$ and $u_iv_i \in E(G)$,

(iii) $u_j = v_i$ and $v_j = u_i$ if $j > i$.

We can also say that if $uv$ is an edge of $S^n_G$, there is
an edge $xy$ of $G$ such that the labels of $u$ and $v$ are: $u = wxyy \ldots y$, $v = wyxx \ldots x$,
where $x,y\in [k]$ and $w \in [k]^\ell$,  $0 \le \ell \le n - 1$.

The generalized Sierpi\'nski graph  $S^n_G$ 
can also be constructed recursively from $k$ copies of  $S^{n-1}_G$ as follows:

- for each $j \in [k]$ add the label $j$ in front of the labels of all vertices in 
$S^{n-1}_G$  and denote the obtained graph by $jS^n_G$,

- for any edge $xy$ of $G$, add an edge between the vertices $xyy\ldots y$ and $yxx \ldots x$ in $S^{n}_G$.

Generalized Sierpi\'nski graphs $S^n_G$ are a natural generalization of {\em Sierpi\'nski graphs with base $p$},  $S^n_p$,
which are generated from complete graphs. In other words, the generalized Sierpi\'nski graph $S^n_{K_p}$ coincides with $S^n_p$. 
Partially motivated by the fact that Sierpi\'nski graphs  belong to a family of subcubic graphs, 
 \cite{BKR} determined bounds on the packing chromatic number of Sierpi\'nski graphs with base 3. 
The exact chromatic number on this class of graph was recently determined (see \cite{ZeVe}).

The packing  colorings of  generalized Sierpi\'nski graphs
as well as of Sierpi\' nski triangle graphs  have been studied by
\cite{BrFe} who determined  the packing chromatic numbers 
of $S^n_G$  for all connected graphs $G$ on 4 vertices
with the exception of two families (generated from $K_4-e$ and the paw graph),  while for  the packing chromatic number of Sierpi\' nski triangle graphs
an upper bound 31 was established.

This paper is organized as follows. 
In the next section, we describe   basic definitions as well as the result which provides the upper bounds on the packing chromatic number of generalized Sierpi\'nski graphs
and Sierpi\'nski triangle graphs.
In Section 3, we report on the  packing chromatic number of  
 generalized  Sierpi\'nski graphs $S^n_G$ where $G$ is a connected graph of order four.
With the presented results, we solve the problem of determining the packing chromatic number for this family of graphs.
 In Section 4, we consider the  packing chromatic number of   
 generalized Sierpi\'nski graphs $S^n_G$, where $G$ is a path or a cycle. These numbers are determined for all 
 paths and cycles with the exception of a cycle of length five where the exact numbers are found till the dimension six 
 while for other dimensions an upper bound is provided. The paper is concluded with results on 
 the packing chromatic number in the family
of Sierpi\'nski--triangle graphs $ST_3^n$. In particular, we show that this number  is bounded from above by 20 
and therefore substantially improve previous results.

\section{Preliminaries}

Let $j \in V(G)$.
If $u \in V(S^{n-1}_G)$, then the corresponding vertex in  $S^{n}_G$ (a ``copy'' of $u$ in $jS^{n-1}_G$) is of the form $ju$. 
Vertices of $S^{n}_G$ of the form $i^n$, $i \in[k]$, are called the {\em extreme vertices}.
Note that only the extreme vertices of any $S^{n-1}_G$ can be end-vertices of edges between distinct copies 
of $S^{n-1}_G$ in $S^{n}_G$. It follows that for $u,v \in V(S^{n-1}_G)$  and $ju, jv \in V(S^{n}_G)$  
we have $d_{S^{n-1}_G}(u,v)=d_{S^{n}_G}(ju,jv)$ (see also \cite{KZ}).

If $ij \in E(G)$, let $^{ij}S^{\ell}_G$ be the graph obtained from $S^{\ell}_G$ by adding the edge between
the extreme vertices  $i^{\ell}$ and $j^{\ell}$. 
We say that $f$ is an {\em extendable packing $c$-coloring of $S^{\ell}_G$} if $f$ is a packing $c$-coloring of 
$^{ij}S^{\ell}_G$  for every $ij \in E(G)$. 

If $H_1$ and $H_2$  are subgraphs of a graph $G$, let $d_G(H_1,H_2)$ denote the distance between $H_1$ and $H_2$, i.e.\ the minimal distance between a vertex of $H_1$ and a vertex of $H_2$ in $G$.

\begin{prop1} \label{expand}
Let $G$ be an undirected graph with vertex set  $[k]$, $i,j \in [k]$, and 
$n > \ell$.
If   a generalized Sierpi\'nski graph $S^{\ell}_G$ admits an extendable packing $c$-coloring
such that  $d_{S^{\ell+1}_G}(iS^{\ell}_G ,j S^{\ell}_G ) > c$ 
for every  $ij \not \in E(G)$, $i\not =j$, 
then $\chi_{\rho}(S^{n}_G) \le c$.
\end{prop1}

\begin{proof}
Let $f$ be an extendable packing $c$-coloring of $S^{\ell}_G$,
and let $f'$ be a $c$-coloring of $S^{\ell +1}_G$ such that $f'$ restricted to $jS^{\ell}_G$ equals $f$ 
for every $j \in E(G)$.

We first show that $f'$  is a packing $c$-coloring of $S^{\ell+1}_G$. 
Let $u=iu', v=jv' \in V(S^{\ell+1}_G)$  where    $u', v' \in V(S^{\ell}_G)$  and $i,j \in  V(G)$.
We have to show that for $f(u)=f(v)=t$, $t \le c$, we have $d_{S^{\ell+1}_G}(u,v) > t$. 

If $i=j$, then $iu'$ and $jv'$ belongs to the same copy of    $S^{\ell}_G$ and since 
$f$ is a packing $c$-coloring of $S^{\ell}_G$, the claim readily follows. 
If $i\not =j$ and $ij \not \in E(G)$, 
then by  $d_{S^{\ell+1}_G}(iS^{\ell}_G ,j S^{\ell}_G ) > c$, we have that
$d_{S^{\ell+1}_G}(iu',jv') > c$. Finally,
if $ij \in E(G)$, 
we claim that $d_{S^{\ell+1}_G}(iu',jv') \ge d_{^{ij}S^{\ell}_G}(u',v') > t$.
Note that $ij^{\ell}$ and  $ji^{\ell}$ are the only vertices of $iS^{\ell}_G$ and $jS^{\ell}_G$, respectively, that  are connected  with an edge.  It follows that
$$d_{S^{\ell+1}_G}(iu',jv') = d_{S^{\ell+1}_G}(iu',ij^{\ell}) + 1 + d_{S^{\ell+1}_G}(ji^{\ell},jv').$$
Since $j^{\ell}$ and  $i^{\ell}$ are adjacent in $^{ij}S^{\ell}_G$, we have
$$d_{^{ij}S^{\ell}_G}(u',v') \le d_{S^{\ell}_G}(u',j^{\ell}) + 1 + d_{S^{\ell}_G}(i^{\ell},v').$$
Equalities $d_{S^{\ell+1}_G}(iu',ij^{\ell}) = d_{S^{\ell}_G}(u',j^{\ell}) $ and
$d_{S^{\ell+1}_G}(ji^{\ell},jv') = d_{S^{\ell}_G}(i^{\ell},v') $ now yield the assertion.
It follows that $f'$  is a packing $c$-coloring of $S^{\ell+1}_G$.

We next show that $f'$  is a packing $c$-coloring of $^{ij}S^{\ell+1}_G$.  Let $i,j \in V(G)$,
$u', v' \in V(S^{\ell}_G)$ and $u=iu', v=jv' \in V(S^{\ell+1}_G)$.
We show  that for $f(u)=f(v)=t$, $t \le c$, we have  
$d_{^{ij}S^{\ell+1}_G}(iu',jv') \ge d_{^{ij}S^{\ell}_G}(u',v') > t$.
By the same argument as above, it suffices to show the claim only for $ij \in E(G)$.
Note that  $^{ij}S^{\ell+1}_G$ is the graph obtained from $S^{\ell+1}_G$ 
 by adding the edge between $i^{\ell+1}$ and $j^{\ell+1}$.
 If a  shortest path  between $iu'$ and $jv'$ does not contain the edge between $i^{\ell+1}$ and $j^{\ell+1}$, 
then $d_{^{ij}S^{\ell+1}_G}(iu',jv') = d_{S^{\ell+1}_G}(iu',jv') $ and by
the discussion above the claim follows.  Otherwise, we have
$$d_{^{ij}S^{\ell+1}_G}(iu',jv') =  d_{S^{\ell+1}_G}(iu',i^{\ell+1}) + 1 + 
d_{S^{\ell+1}_G}(j^{\ell+1},jv').$$
Since $j^{\ell}$ and  $i^{\ell}$ are adjacent in $^{ij}S^{\ell}_G$, we have
$$d_{^{ij}S^{\ell}_G}(u',v') \le d_{S^{\ell}_G}(u',j^{\ell}) + 1 + d_{S^{\ell}_G}(i^{\ell},v').$$
Equalities $d_{S^{\ell+1}_G}(iu',i^{\ell+1}) = d_{S^{\ell}_G}(u',i^{\ell}) $ and
$d_{S^{\ell+1}_G}(j^{\ell+1},jv') = d_{S^{\ell}_G}(j^{\ell},v') $ now yield the assertion.

Since $f'$  is a packing $c$-coloring of $^{ij}S^{\ell+1}_G$ for every 
$ij \in E(G)$, we showed that $f'$ is an extendable packing $c$-coloring of $S^{\ell+1}_G$.
In other words, we showed that a packing $c$-coloring of $S^{n}_G$ can be obtained by 
using $f$ for every copy  of $S^{\ell}_G$ in $S^{n}_G$.  
This assertion  completes the proof.
\end{proof}

The questions of packing coloring for various  finite and infinite graphs  have been reduced to SAT problems already 
in \cite{zehvesel}, \cite{inglezi} and \cite{KoVe}. 
In this paper,  we applied this approach for  searching packing $c$-colorings and extendable packing $c$-colorings of 
generalized  Sierpi\'nski  graphs as well as  Sierpi\'nski triangle graphs.  
In particular, we applied propositional formulas presented in \cite{KoVe} which
transform an instance of a packing coloring problem into a propositional satisfiability
test (SAT). 
We used the SAT-solvers Glucose Syrup 4.1 (see \cite{Glucose}) and Cryptominisat 5 (see \cite{Cryptominisat})  
in order to find the solutions of the derived 
propositional formulas.

Some of the constructed colorings that provide improved upper bounds  are  not present in the paper. 
Interested readers are invited to visit the website 
\url{https://omr.fnm.um.si/wp-content/uploads/2017/06/SierpinskiP.pdf} 
where all of the obtained colorings are given.



\section{Generalized  Sierpi\'nski  graphs with base graphs on 4 vertices}

There are altogether six generalized  Sierpi\'nski  graphs with base graphs on 4 vertices. 
As mentioned previously, 
the packing chromatic numbers of $S^n_G$  
for all connected graphs $G$ on 4 vertices have been determined in \cite{BrFe}
with the exception
of $K_4 -  e$ (the graph obtained by removing an edge from $K_4$) and the {\em paw graph} 
(the graph obtained by joining  one vertex of $K_3$ to a  $K_1$).

\begin{figure}[hbt] 
\begin{center}

\includegraphics[width=110mm]{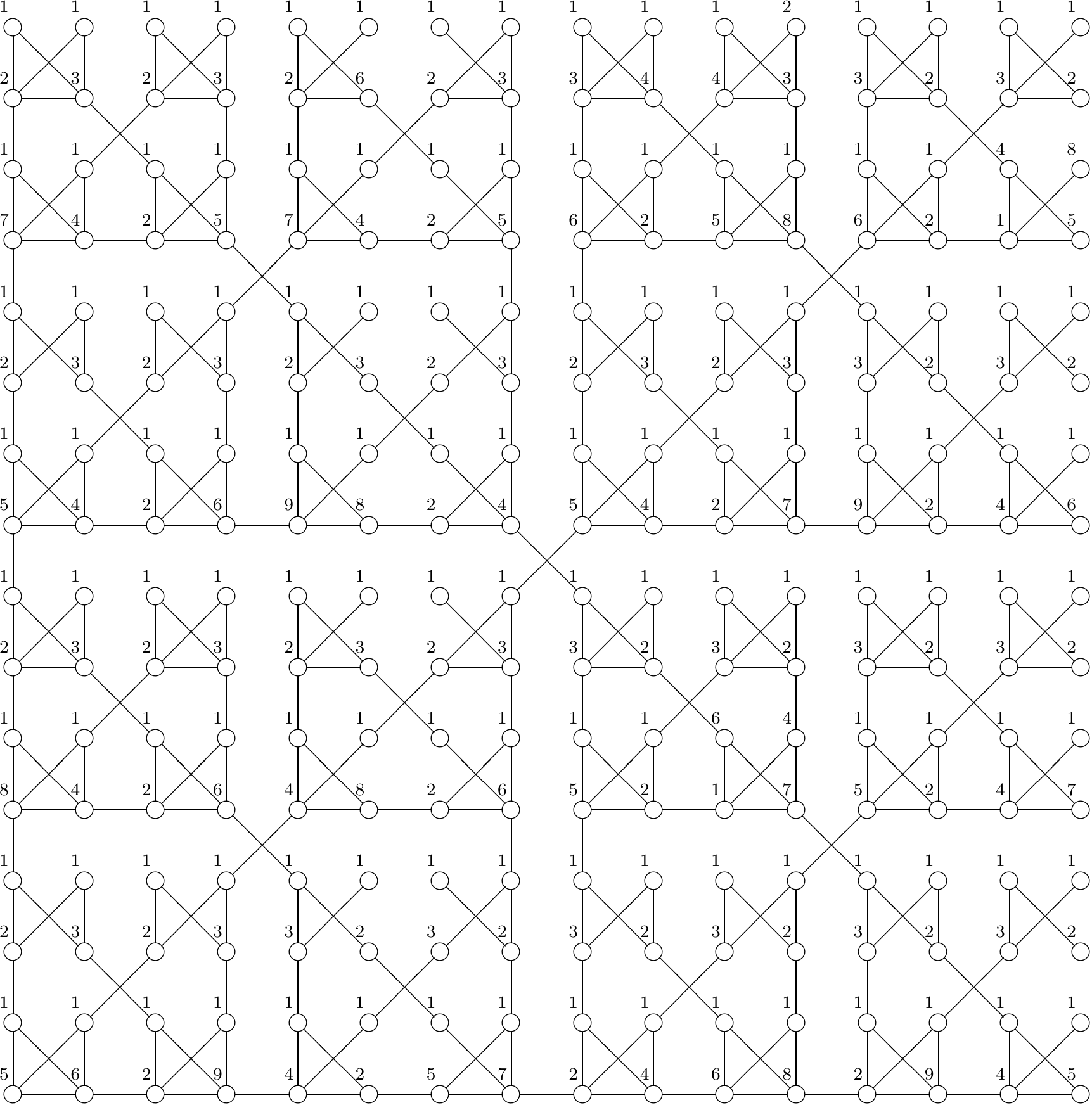}  

\caption{A packing 9-coloring of $S_{K_4-e}^4$   } \label{9k4-e4}
\end{center}
\end{figure}

 In this section, we establish the packing chromatic numbers for these two families of graphs.

\begin{thm1}
If $S_{K_4-e}^n$ is  
the generalized Sierpi\'nski  graph of $K_4 -  e$ of dimension $n$, then
\begin{displaymath}
 \chi_\rho(S_{K_4-e}^n)  =  \left\{ \begin{array}{ll}
 3,  &  \;  n=1 \\ 
 6,  &  \;  n=2 \\ 
 8,  &  \;  n=3 \\ 
 9,  &  \;  n=4 \\ 
 10,  &  \;  n\ge 5 \\ 
 \end{array} \right.
\end{displaymath}
\end{thm1}

\begin{figure}[hbt] 
\begin{center}

\includegraphics[width=131mm]{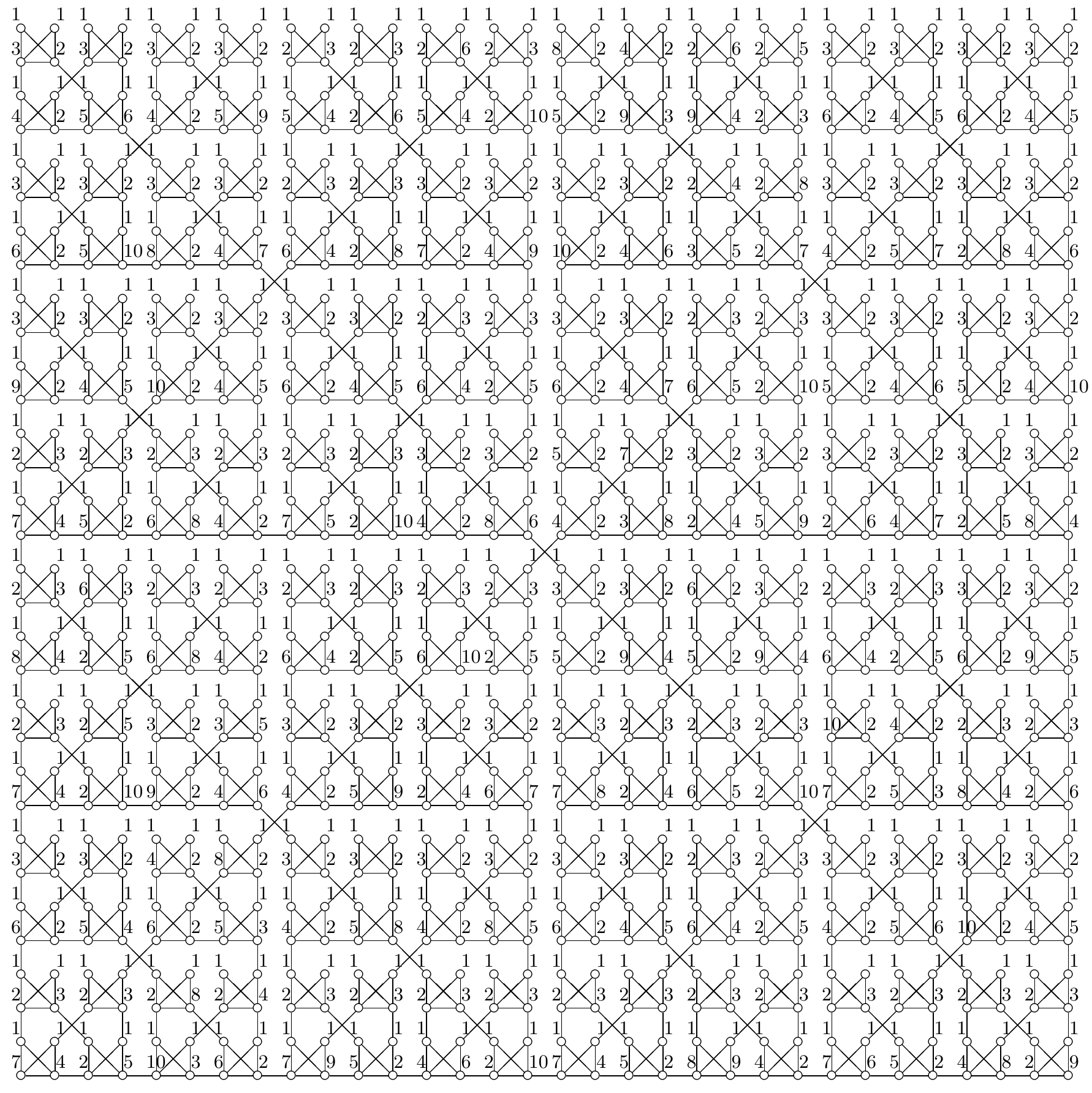}

\caption{An extendable packing 10-coloring of $S_{K_4-e}^5$   }   \label{10k4-e5}
\end{center}
\end{figure}

\begin{proof}
For $n \le 3$, the result is presented in \cite{BrFe}. 
Since we showed with SAT solver that neither a packing 8-coloring of $S_{K_4-e}^4$ 
nor a packing 9-coloring of $S_{K_4-e}^5$ 
can be obtained, we have the lower bound for $n=4$ as well as for $n\ge 5$. 

The upper bound for $n=4$ was obtained 
by a packing 9-coloring of $S_{K_4-e}^4$  presented in Fig. \ref{9k4-e4}, 
while for $n \ge 5$ we obtained an extendable packing 10-coloring of $S_{K_4-e}^5$    presented in Fig. \ref{10k4-e5}.
\end{proof}

\begin{figure}
  \centering
\begin{subfigure}

  \includegraphics[width=5.5cm]{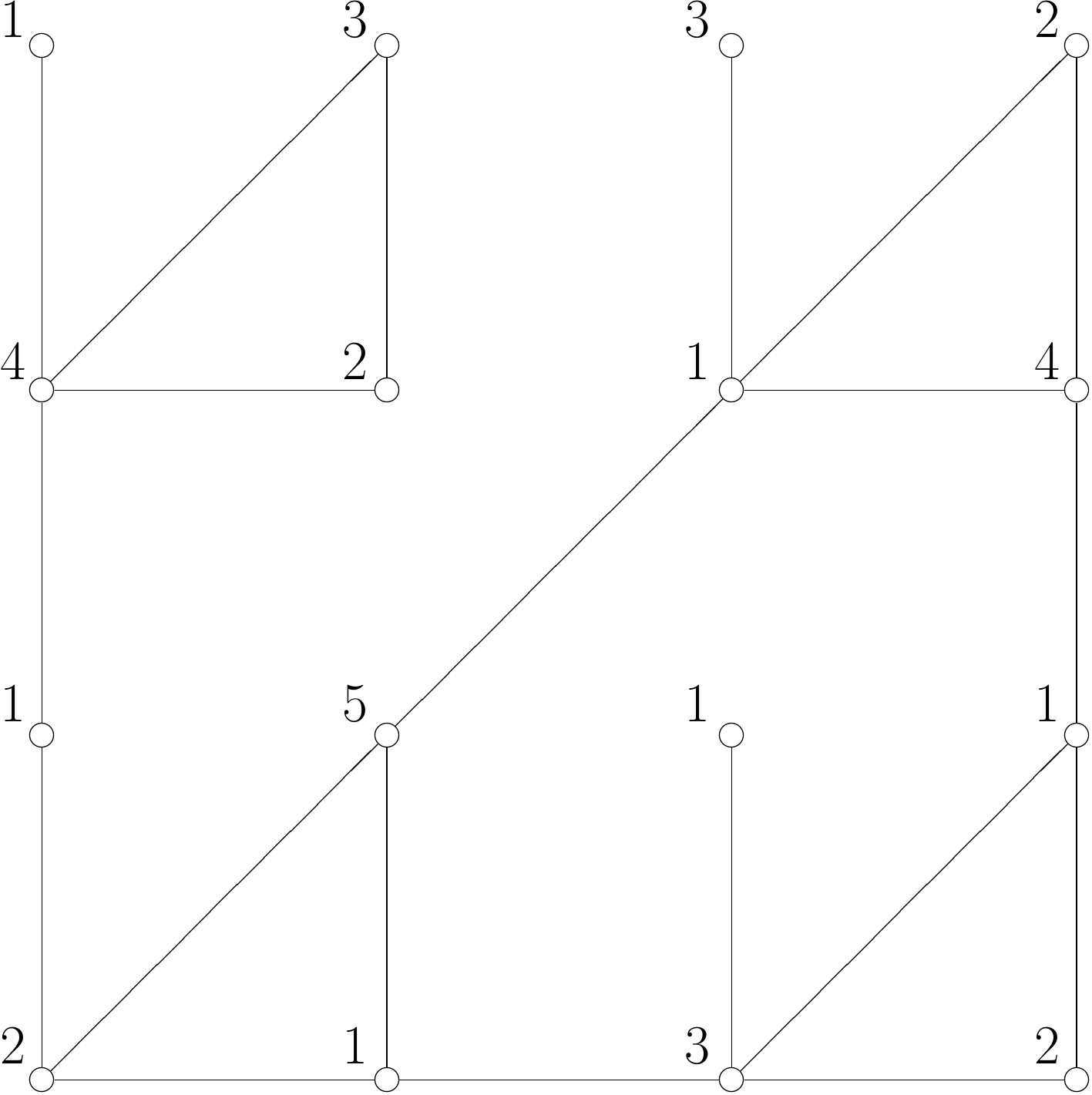}
\end{subfigure}%
\hspace{10mm}
\begin{subfigure}

  \includegraphics[width=5.5cm]{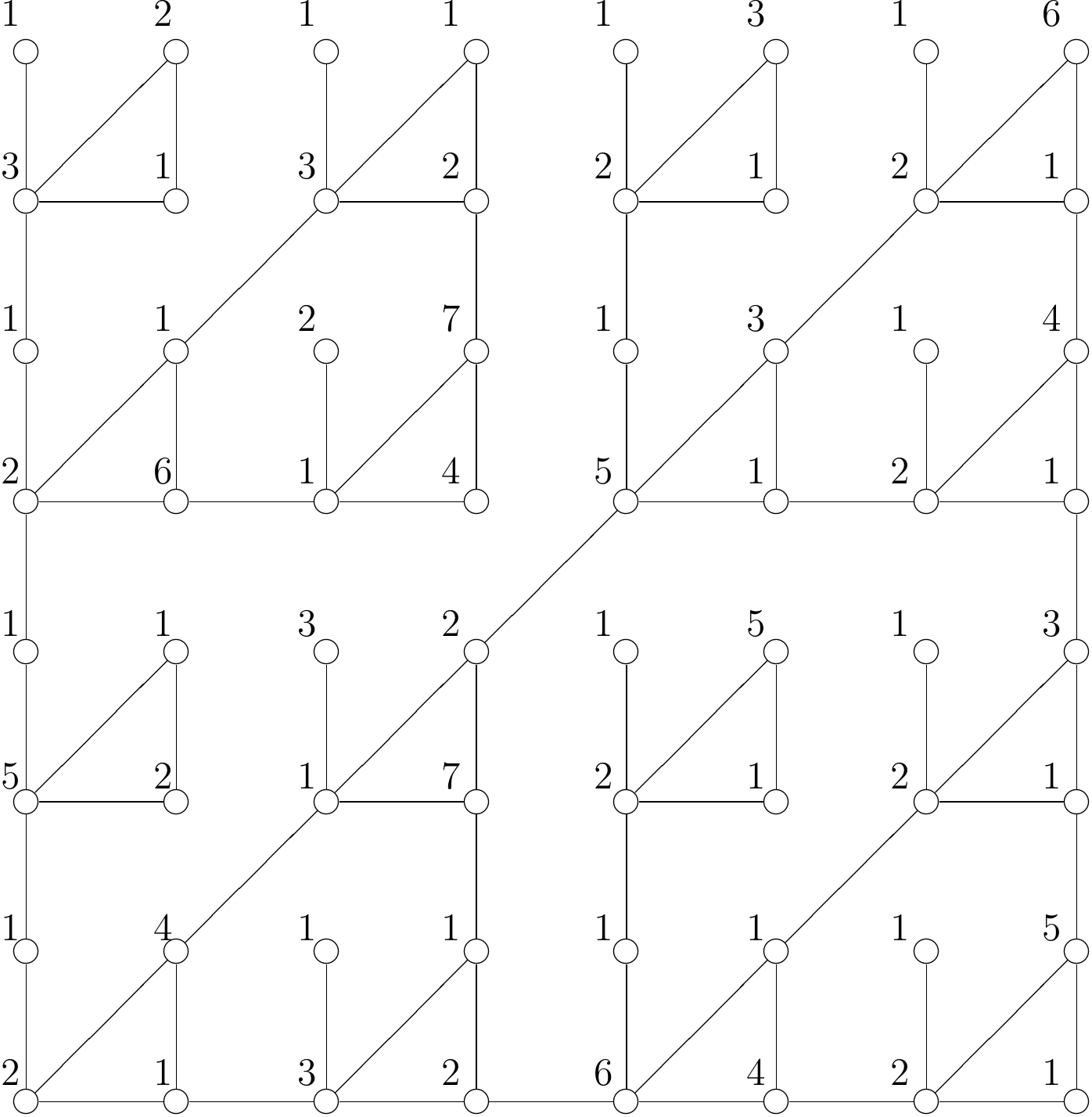}
\end{subfigure}
  \caption{A packing 5-coloring of $S_{paw}^2$ (left) and  packing 7-coloring of $S_{paw}^3$ (right)} \label{7Paw3}
\end{figure}


\begin{thm1}
If $S_{paw}^n$ is  
the generalized Sierpi\'nski  graph of the paw of dimension $n$, then
\begin{displaymath}
 \chi_\rho(S_{paw}^n)  =  \left\{ \begin{array}{ll}
 3,  &  \;  n=1 \\ 
 5,  &  \;  n=2 \\ 
 7,  &  \;  n=3 \\ 
 8,  &  \;  n \ge 4 \\ 
 \end{array} \right.
\end{displaymath}
\end{thm1}

\begin{figure}[hbt] 
\begin{center}

\includegraphics[width=120mm]{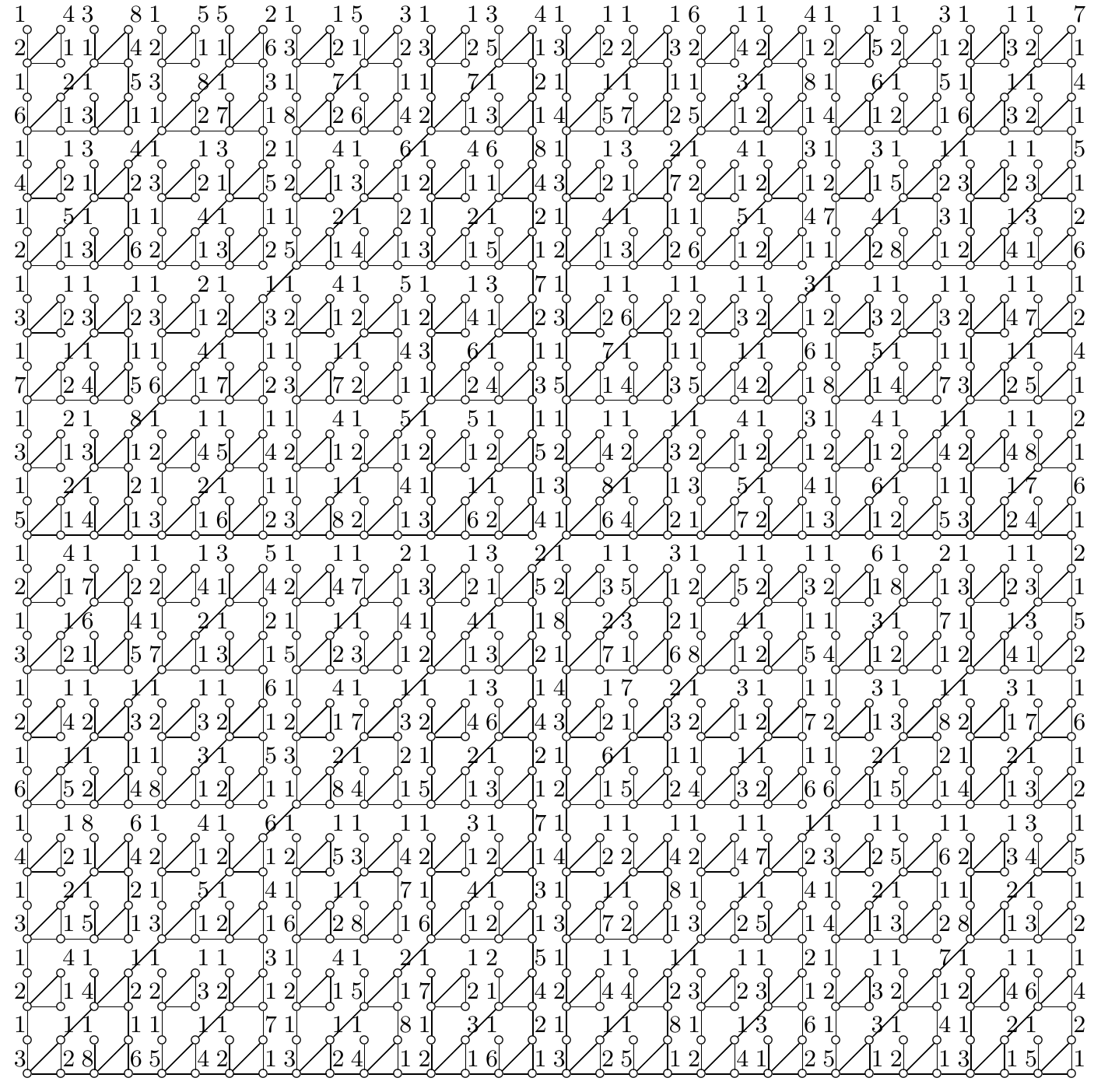} 

\caption{An extendable packing 8-coloring of $S_{paw}^5$   } \label{8Paw5}
\end{center}
\end{figure}

\begin{proof}
For $n = 1$, this is the result presented in \cite{BrFe}, as well as the lower bounds for  $n=2$ and $n=3$. 

Since we showed with our SAT solver that a packing 7-coloring of $S_{paw}^4$ 
cannot be obtained, we have the lower bound for $n\ge 4$. 

The upper bounds were obtained 
by the following constructions: a packing 5-coloring of $S_{paw}^2$ 
and a packing 7-coloring of $S_{paw}^3$   presented in Fig.  \ref{7Paw3}, 
while for $n \ge 4$ we obtained an extendable packing 8-coloring of $S_{paw}^5$    presented in Fig. \ref{8Paw5}.

\end{proof}

\section{Generalized  Sierpi\'nski  graphs with base graphs paths and cycles}

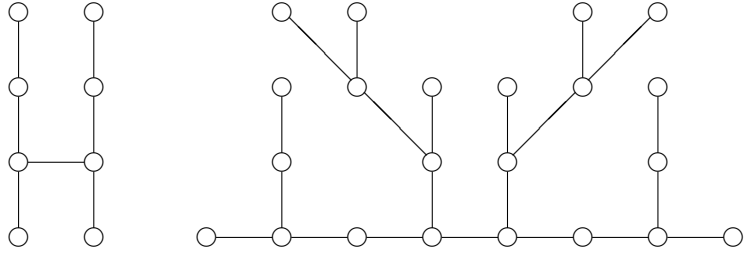
\begin{figure}[hbt] 
\begin{center}  \unitlength = 0.25mm

\begin{picture}(400,150) (0,0)

\put(0,0){\circle{10}}
\put(0,40){\circle{10}}
\put(0,80){\circle{10}}
\put(0,120){\circle{10}}
\put(40,0){\circle{10}}
\put(40,40){\circle{10}}
\put(40,80){\circle{10}}
\put(40,120){\circle{10}}

\put(0,5){\line(0,1){30}}
\put(0,45){\line(0,1){30}}
\put(0,85){\line(0,1){30}}
\put(40,5){\line(0,1){30}}
\put(40,45){\line(0,1){30}}
\put(40,85){\line(0,1){30}}

\put(5,40){\line(1,0){30}}

\put(100,0){\circle{10}}
\put(105,0){\line(1,0){30}}
\put(140,0){\circle{10}}
\put(145,0){\line(1,0){30}}
\put(180,0){\circle{10}}
\put(185,0){\line(1,0){30}}
\put(220,0){\circle{10}}
\put(225,0){\line(1,0){30}}
\put(260,0){\circle{10}}
\put(265,0){\line(1,0){30}}
\put(300,0){\circle{10}}
\put(305,0){\line(1,0){30}}
\put(340,0){\circle{10}}
\put(345,0){\line(1,0){30}}
\put(380,0){\circle{10}}

\put(140,40){\circle{10}}
\put(140,5){\line(0,1){30}}
\put(220,40){\circle{10}}
\put(220,5){\line(0,1){30}}
\put(260,40){\circle{10}}
\put(260,5){\line(0,1){30}}
\put(340,40){\circle{10}}
\put(340,5){\line(0,1){30}}

\put(140,80){\circle{10}}
\put(180,80){\circle{10}}
\put(140,45){\line(0,1){30}}
\put(220,80){\circle{10}}
\put(220,45){\line(0,1){30}}
\put(217,44){\line(-1,1){33}}
\put(260,80){\circle{10}}
\put(263,44){\line(1,1){33}}
\put(300,80){\circle{10}}
\put(260,45){\line(0,1){30}}
\put(340,80){\circle{10}}
\put(340,45){\line(0,1){30}}

\put(140,120){\circle{10}}
\put(180,120){\circle{10}}
\put(180,85){\line(0,1){30}}
\put(177,84){\line(-1,1){33}}

\put(340,120){\circle{10}}
\put(300,120){\circle{10}}
\put(300,85){\line(0,1){30}}
\put(303,84){\line(1,1){33}}

\end{picture}

\caption{Graph  $H'$ (left) and  $H$ (right)  } \label{hh}
\end{center}
\end{figure}

\begin{prop1} \label{spodnjaPk}
Let $n$ and $k$ be integers with $k\ge 4$ and $n\ge 3$.
If $S_{P_k}^n$ is  
the generalized Sierpi\'nski  graph of  $P_k$ of dimension $n$, then

  (i) $\chi_\rho(S_{P_k}^2)  \ge 4$,

  (ii) $\chi_\rho(S_{P_k}^n)  \ge 5$.
\end{prop1}

\begin{proof}
(i) Note that $A=\{11, 12, 13, 14, 20, 21, 22, 23  \}$ is a subset of $V(S_{P_k}^2)$. 
It can be seen that the graph induced by $A$ is isomorphic to the graph $H'$ depicted
in the left-hand side of Fig \ref{hh}.
Since it shown in \cite[Lemma 4]{BrFe} that $\chi_\rho(H') \ge 4$, this  case is settled. 

(ii) The set $B=\{111, 112, 113, 121, 120, 122, 123, 132, 131, 133, 134,
211, 210, 201, 200, 202, 203,$ $ 212, 213, 221, 222 \}$ is a subset of $V(S_{P_k}^3)$. 
Moreover,  the graph induced by $B$ is isomorphic to the graph $H$
depicted in the right-hand side of Fig \ref{hh}.
Since it shown in \cite[Lemma 4]{BrFe} that 
$\chi_\rho(H) \ge 5$, the assertion follows. 
\end{proof}

\begin{prop1} \label{Pk2} 
Let $k \ge 3$.
If $S_{P_k}^2$ is the  generalized Sierpi\'nski  graph of  $P_k$ of dimension $2$, then 

\begin{displaymath}
 \chi_\rho(S_{P_k}^2)  =  \left\{ \begin{array}{ll}
 3,  &  \;  k=3 \\ 
 4,  &  \;  k \ge 4 \\ 
 \end{array} \right.
\end{displaymath}

\end{prop1}

\begin{proof}
Since $S_{P_3}^2$ induces a path with seven vertices, we have $\chi_\rho(S_{P_k}^2)  \ge 3$.
Let $f$ be a function from $V(S_{P_3}^2)$ to the set of integers, defined by $f(00)=f(02)=f(10)=f(12)=f(20)=f(22)=1$,
$f(01)=f(21)=2$ and $f(11)=3$. Obviously,  $f$ is a packing 3-coloring of $S_{P_3}^2$ and this case is settled.

For $k\ge 4$, the lower bound is given by Proposition \ref{spodnjaPk}. 
In order to prove the upper bound, note that 
by the definition of the  generalized Sierpi\'nski  graph,  $S_{P_k}^2$ is composed of $k$ copies of $P_k$ such that 
an edge joins a vertex $x$ in $iP_k$ with a vertex $y$ in  $jP_k$  if and only if $i=j\pm 1$. 
The edge $xy$ is called a {\em cross edge}. 
Note that  $x$ and $y$ are of the form $ij$ and $ji$, respectively. 



Let $f$ be a 4-coloring of  $S_{P_k}^2$ defined as follows:

\begin{displaymath}
 f(ij)  =  \left\{ \begin{array}{ll}
 1,  &  \;  j \equiv  0 \;({\rm mod} \; 2 )     \\ 
 2,  &  \;  j \equiv  1\; ({\rm mod} \; 4 )  \; {\rm and} \; i \not = j  \\ 
 3,  &  \;  j \equiv  3\; ({\rm mod} \; 4 ) \; {\rm and} \; i \not = j  \\ 
 4,  &  \; j \equiv 1  \;({\rm mod} \; 2 ) \; {\rm and} \; i  = j \\ 
 \end{array} \right.
\end{displaymath}

If $i$ is even, then $f$ restricted to $iP_k$ is the sequence $1,2,1,3,1,2,1,3,\ldots$, while  
for $i$  odd, we obtain the same sequence with the exception that it 
admits exactly one vertex with color 4 at the position $i$. 
It is straightforward to see, that both sequences imply  a packing 4-coloring of $P_k$.  
In order to show that $f$ is a packing 4-coloring of $S_{P_k}^2$, 
we have to consider vertices which are ``close'' to a cross edge $xy$.
If $i\in [k-1]$, we may assume without loss of generality that $x=ij=i(i+1)$  and $y=ji=(i+1)i$.
Let  $j\ge 1$, $j\ge 2$, $j\le k-2$ and $j\le k-3$ for $x^- := i(j-1)$, $x^{--} := i(j-2)$, 
$x^+ := i(j+1)$ and $x^{++} := i(j+2)$, respectively.
Analogously, we define   $y^{--}, y^{-},y^{+}$ and $y^{++}$. 
Since  either $f(x)=1$ or $f(y)=1$, for all vertices $u$ and   $v$, $u\not = v$, with  $f(u)=f(v)=1$
we have $d(u,v)\ge 2$. Moreover, by definition of $f$, for $f(u)=f(v)=4$ we have $d(u,v)\ge 6$. 
It follows that we have to study only  vertices with color 2 or 3. 

By definition of $f$, for $x=ij$ and $y=ji$ we have 

\begin{displaymath}
 f(x)  =  \left\{ \begin{array}{ll}
 1,  &  \;  j \equiv  0 \;({\rm mod} \; 2 )         \\ 
 2,  &  \;  j \equiv  1 \;({\rm mod} \; 4 )         \\ 
  3,  &  \;  j \equiv  3 \;({\rm mod} \; 4 )   \\ 
 \end{array} \right.
\end{displaymath}
  
and 

\begin{displaymath}
 f(y)  =  \left\{ \begin{array}{ll}
 1,  &  \;  j \equiv  1 \;({\rm mod} \; 2 )     \\ 
 2,  &  \;  j \equiv  2 \;({\rm mod} \; 4 )    \\ 
 3,  &  \;  j \equiv  0 \;({\rm mod} \; 4 )   \\ 
 \end{array} \right.
\end{displaymath}
  
Thus, we consider  the following cases. 

(i) If $f(x)=1$, then  either $f(y)=2$ or  $f(y)=3$. We have that $f(x^-)=4$ and $f(x^{--})=f(x^{++})=1$.
If  $f(y)=2$ (resp.  $f(y)=3$), then  $f(x^{+})=3$ (resp. $f(x^{+})=2$). 
The vertices closest to  $x^{+}$ in $(i+1)P_k$ with color 3 (resp. 2) are $y^{--}$ and $y^{++}$. 
Since both are at distance four from  $x^{+}$, this case is settled.

(ii) If $f(y)=1$, then either  $f(x)=3$ or $f(x)=2$. 
 We have that $f(x^+)=f(x^-)=1$.
If $f(x)=3$  (resp. $f(x)=2$), then   $f(x^{--})=f(x^{++})=2$ (resp. $f(x^{--})=f(x^{++})=3$). 
Since $d(x^{--},y)=d(x^{++},y)=3$, we have to consider only vertices close to $x$. 
Note that $f(y^-)=4$, thus, vertices in $(i+1)P_k$ of color 3 (resp. 2) are clearly at distance at least four from $y$ and 
the assertion follows. 
\end{proof}

\begin{thm1} \label{Pkn}
Let $k$ and $n$ be integers.
If $S_{P_k}^n$ is the  generalized Sierpi\'nski  graph of  $P_k$ of dimension $n$, then 
\begin{displaymath}
 \chi_\rho(S_{P_k}^n)  =  \left\{ \begin{array}{ll}
 3,  &  \;  k=3 \; {\rm and} \;  n \ge 2\\ 
 4,  &  \;  k \ge 4 \; {\rm and} \;  n = 2 \\
 5,  &  \;  k \ge 4 \; {\rm and} \;  n \ge 3 \\
  
 \end{array} \right.
\end{displaymath}

\end{thm1} 

\begin{proof}
For $n=2$ this is Proposition \ref{Pk2}, while for $n\ge 3$ the lower bound is given by Proposition \ref{spodnjaPk}. 
If $k=3$ and $n\ge 3$ note that the packing 3-coloring of $S_{P_3}^2$ presented  in the proof of Proposition  \ref{Pk2}  is 
an extendable  packing 3-coloring of this graph.

Let $n \ge 3$, $s,i,j \in [k]$ and $w \in [k]^{n-3}$. 
We define  a 5-labeling $f$  of $S_{P_k}^n$ as follows:

\begin{displaymath}
 f(wsij)  =  \left\{ \begin{array}{ll}
 1,  &  \; j \equiv  0 \;({\rm mod} \; 2 )     \\ 
 2,  &  \;  j \equiv  1 \;({\rm mod} \; 4 )  \; {\rm and} \; i \not = j  \\ 
 3,  &  \;  j \equiv  3 \;({\rm mod} \; 4 )   \; {\rm and} \; i \not = j  \\ 
 4,  &  \;  j \equiv  1 \;({\rm mod} \; 2 )  \; {\rm and} \; i  = j   \; {\rm and} \; s \equiv  0 \;({\rm mod} \; 2 )  \\ 
 5,  &  \; j \equiv  1 \;({\rm mod} \; 2 )  \; {\rm and} \; i  = j   \; {\rm and} \; s \equiv  1 \;({\rm mod} \; 2 )  \\ 
 \end{array} \right.
\end{displaymath}

From  the proof of Proposition \ref{Pk2} it follows that a restriction of $f$ to a copy of $S_{P_k}^2$ in $S_{P_k}^n$  is 
either a packing 4-coloring or a packing 5-coloring of $S_{P_k}^2$. 

In order to see that  $f$ is  
a packing 5-coloring of $S_{P_k}^n$, note that an edge of $S_{P_k}^n$
that does not belong to a copy of $S_{P_k}^2$ connects  vertices $u$ and $v$
of the form $u=wij^{n-\ell-1}$, $v=wji^{n-\ell-1}$, $i,j \in [k]$, $0\le \ell \le n-3$, $w \in [k]^{\ell}$ and $|i - j| =1$. 

Assume that $f$ admits vertices $\alpha \in V(zjS_{P_k}^2)$ and   $\beta \in V(wiS_{P_k}^2)$ such that $f(\alpha)=f(\beta)=\xi$ and 
$d(\alpha,\beta) \le \xi$.  We may assume without loss of generality  that $ i \equiv  0 \;({\rm mod} \; 2 ) $.
Thus, by definition of $f$, we have $f(u)=1$ and 
$f(v) \in \{4,5\}$. Note that $v$ is of degree three, say $N(v)= \{u, u', u''\}$. Moreover, all vertices in $N(v)$   have color 1. 
It follows that $\xi \not =1$. 
If $f(v)=4$ (resp. $f(v)=5$),  
then the copy of $S_{P_k}^2$  that contains $u$ does not admit a vertex with color 4 (resp. 5).
It follows that $\xi \not \in \{4,5\}$. 
Finally, since  the vertices of $N(u)\setminus \{v \} $ are at distance four from the vertices of 
$(N(u') \cup N(u'')) \setminus \{ v \} $, we have  $\xi \not \in \{2,3\}$ and we obtain a contradiction.
Since we showed that $f$ is a  packing 5-coloring of $S_{P_k}^n$, the proof is complete.
\end{proof}

\begin{figure}[hbt] 
\begin{center}
\includegraphics[width=100mm]{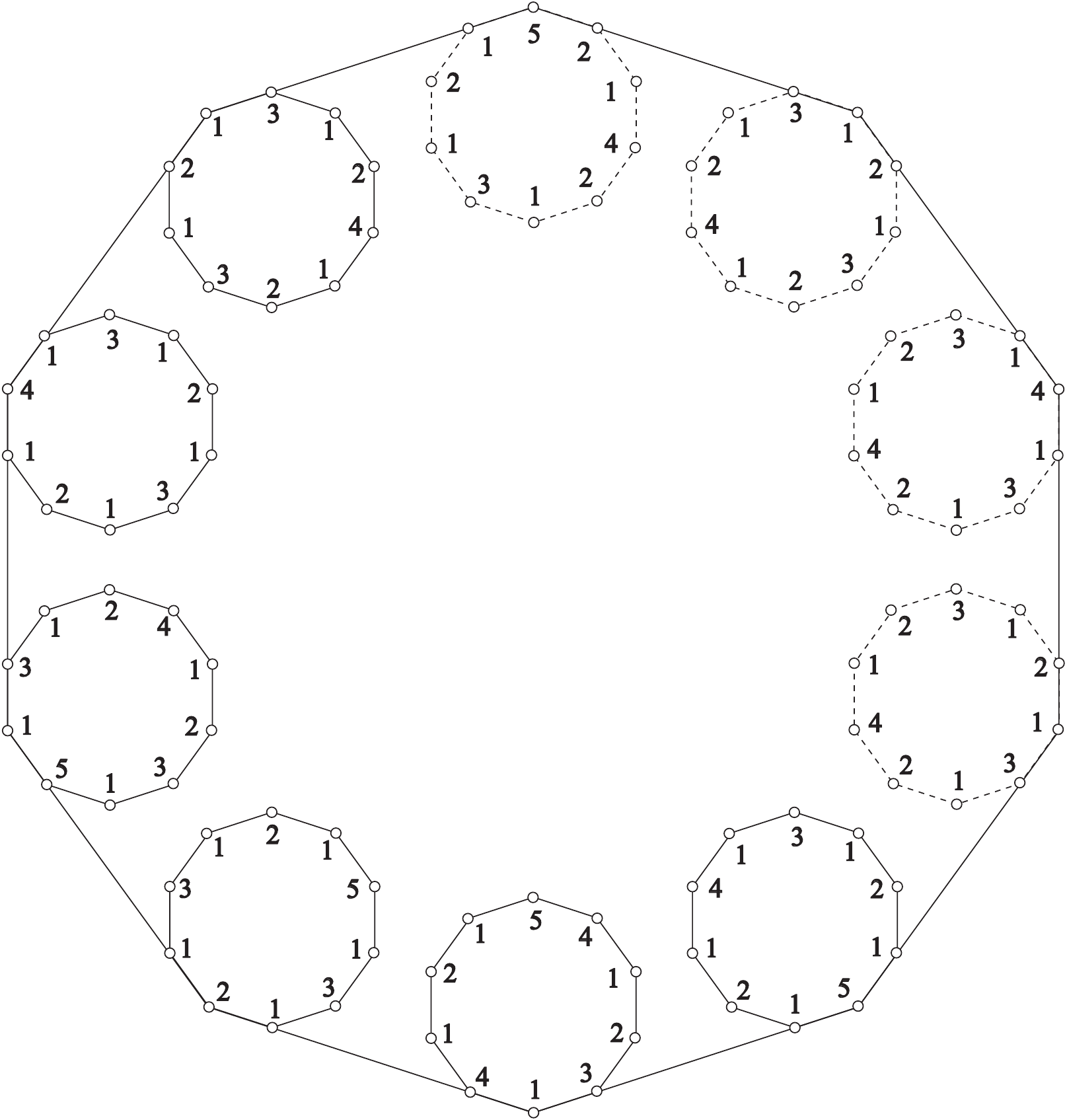} 
\caption{An extendable packing 5-coloring of $S_{C_{10}}^2$  } \label{5C10}
\end{center}
\end{figure}

\begin{figure}[hbt] 
\begin{center}
\includegraphics[width=100mm]{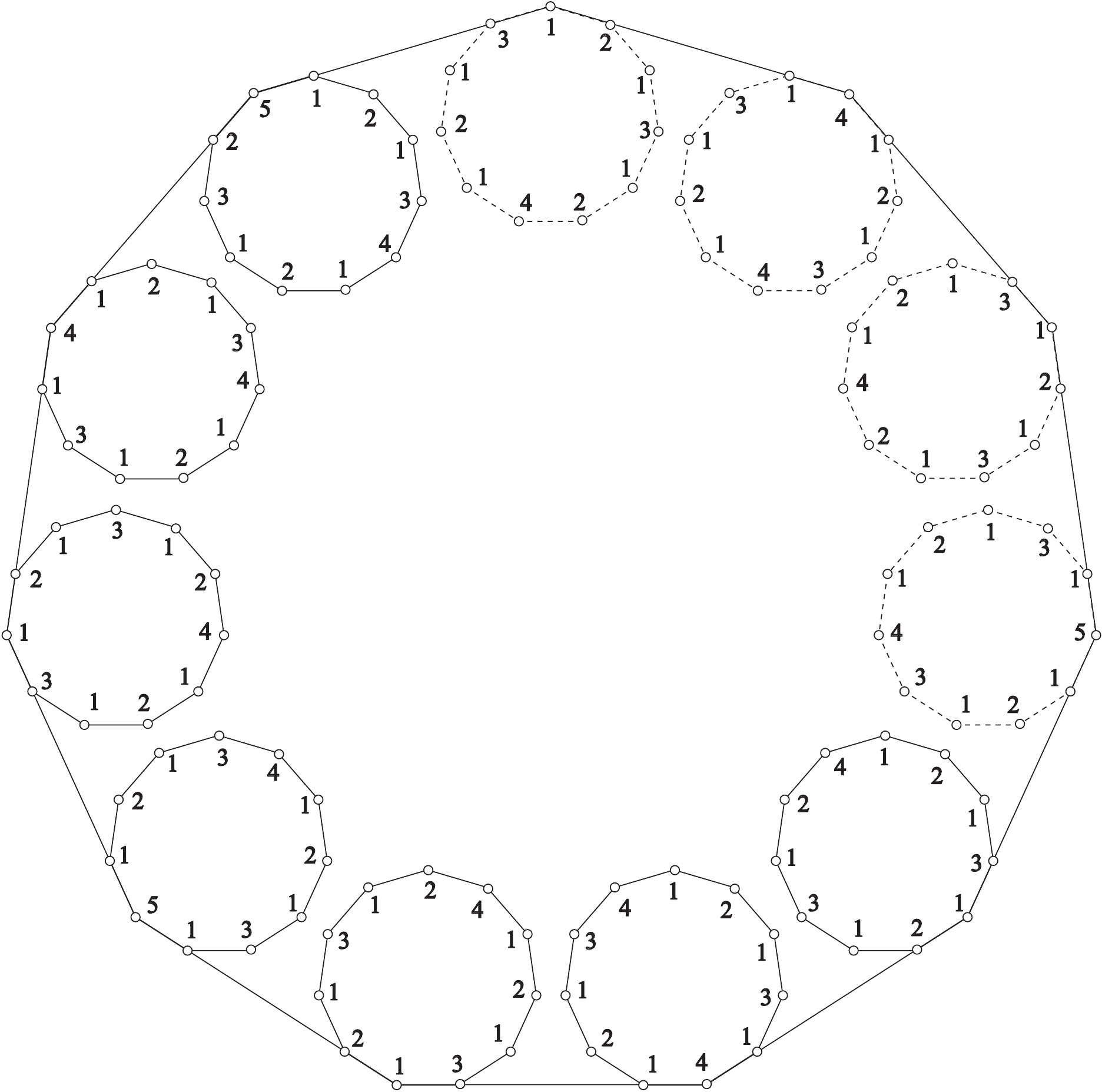} 
\caption{An extendable packing 5-coloring of $S_{C_{11}}^2$ } \label{5C11}
\end{center}
\end{figure}

\begin{figure}[hbt] 
\begin{center}
\includegraphics[width=100mm]{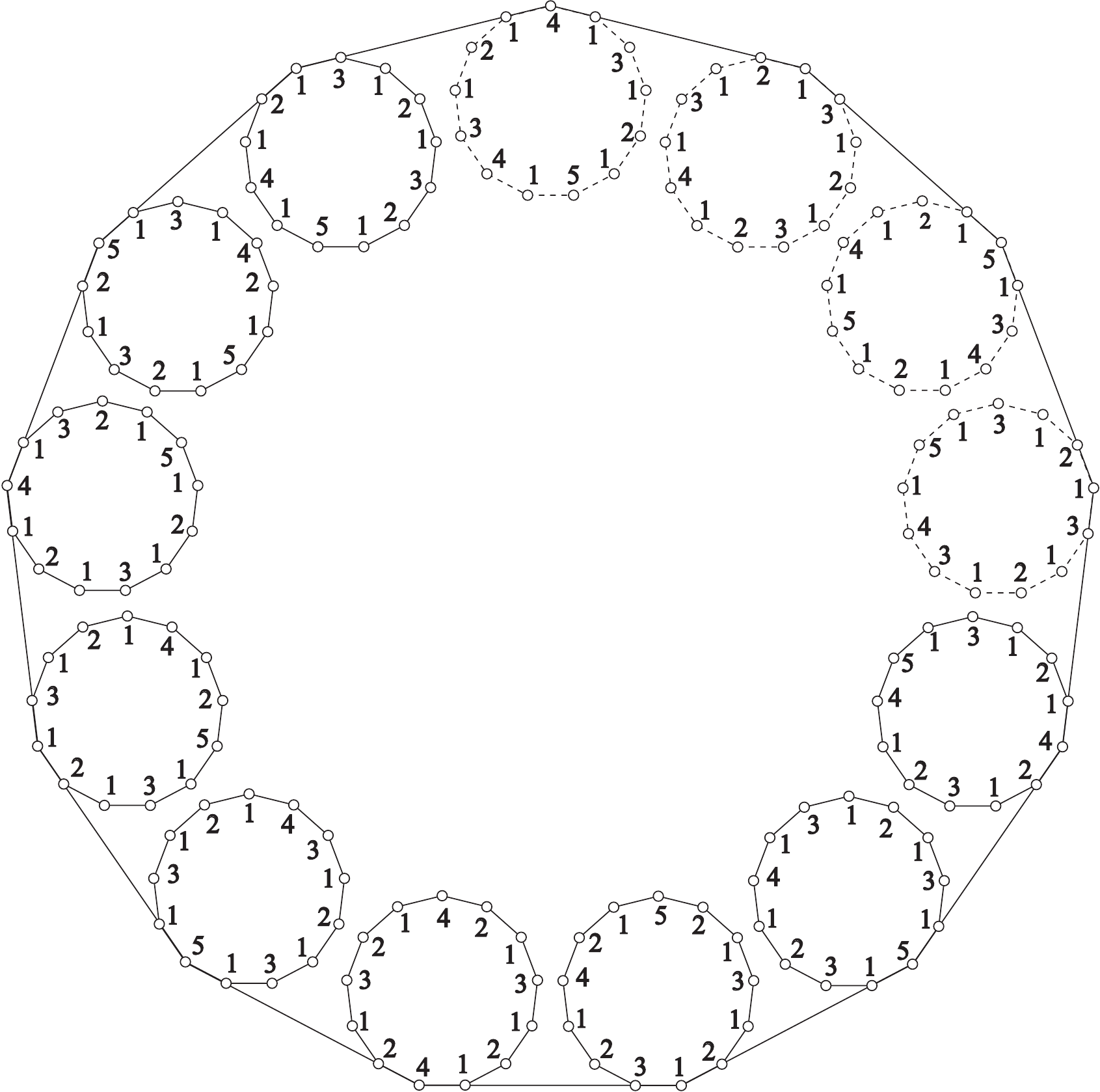} 
\caption{An extendable packing 5-coloring of $S_{C_{13}}^2$  } \label{5C13}
\end{center}
\end{figure}


Let the $n \ge 2$, $i,j \in [k]$ and $w \in [k]^{n-2}$. 
Note that $S_{C_k}^n$ is composed of $k^{n-2}$ copies of $S_{C_k}^2$. 
We will say that vertices of the  form $wij$ belong to a copy of $S_{C_k}^2$ denoted as $wS_{C_k}^2$.
Moreover, $wS_{C_k}^2$  is composed of $k$ copies of $C_k$ and a copy that admits the vertex $wij$ will be denoted as $iC_k$.

In the proofs of  below propositions,  we  define  a packing coloring of $wS_{C_k}^2$ for every $w \in [k]^{n-2}$ such that the coloring of every $iC_k$ is given either explicitly or
as  a sequence of colors  for vertices $wii, wi(i+1), wi(i+2), \ldots, wi(i+k-1)$ (addition modulo $k$). 
 Let also
$\overline{wxyz}$ stand for the sequence of colors $wxyz$ which  can be repeated  as needed 
in order to color all vertices of $iC_k$.

\begin{prop1} \label{uC2}
Let  $k \ge 10$  and $k$ is even or  $k \equiv 0\;(mod \;4)$.  
If $S_{C_k}^2$ is the  generalized Sierpi\'nski  graph of  $C_k$ of dimension $2$, then 
$ \chi_\rho(S_{C_k}^2)  =  4$
\end{prop1}

\begin{proof}
Note that $S_{C_k}^2$ is composed of $k$ copies of $C_k$ such that for $i,j \in [k]$  
an edge joins a vertex $x$ in $iC_k$ with a vertex $y$ in  $jC_k$  if and only if $i=j + 1$ (addition modulo $k$). 
Moreover, $x$ and $y$ are of the form $ij$ and $ji$, respectively. 

If $k \equiv 0 \; ($mod$ \;4) $, then
$f$ is a 4-coloring of $S_{C_k}^2$ defined by:
\begin{displaymath}
 f(ij)  =  \left\{ \begin{array}{ll}
 1,  &  j  \equiv 0 \; ({\rm mod} \; 2 )     \\ 
 2,  &   j  \equiv 1 \; ( {\rm mod} \; 4 )   \; {\rm and} \; i \not = j  \\ 
 3,  &    j \equiv 3 \; ({\rm mod} \; 4 )   \; {\rm and} \; i \not = j  \\ 
 4,  &    j \equiv 1\;   ({\rm mod} \; 2 )   \; {\rm and} \; i  = j \\ 
 \end{array} \right.
\end{displaymath}

For $k\ge 10$ and $k \equiv 2 \; ($mod$ \;4) $, we define  a 4-coloring of $S_{C_k}^2$ as follows.
Let $i,j \in [k]$ and  let $ij$ be a vertex of $S_{C_k}^2$.
We label $iC_k$ such that we start at the vertex $ii$. If $i$ is even  (resp. $i$ is odd), 
 we use the sequence $121314\overline{1213}$ (resp. $4121314213\overline{1213}$).

We can see analogously as in the proof of Proposition \ref{Pk2} that above colorings are packing 4-coloring of $S_{C_k}^2$ for
$k\equiv 0 \; ( mod \;4)$ and  $k\equiv 2 \; ( mod \;4) $, respectively.

\end{proof}

\begin{prop1} \label{uCk}
Let $n\ge 2$ and  $k \not\in \{5, 6, 7\}$. 
If  $S_{C_k}^n$ is the  generalized Sierpi\'nski  graph of  $C_k$ of dimension $n$, then 
$ \chi_\rho(S_{C_k}^n)  \le  5$.
\end{prop1}

\begin{proof}

Let $k=4t$, $s,i,j \in [k]$ and $w \in [k]^{n-3}$. 
We define  a 5-coloring $f$ of $S_{C_{4t}}^n$ as follows:

\begin{displaymath}
 f(wsij)  =  \left\{ \begin{array}{ll}
 1,  &  j \equiv 0  \;  ({\rm mod} \; 2 )    \\ 
 2,  &  j \equiv 1  \;  ({\rm mod} \; 4 ) \; {\rm and} \; i \not = j  \\ 
 3,  &  j \equiv 3  \;  ({\rm mod} \; 4 )  \; {\rm and} \; i \not = j  \\ 
 4,  &  j \equiv 1  \;  ({\rm mod} \; 2 )\; {\rm and} \; i  = j   \; {\rm and} \;  s \equiv  0 \;({\rm mod} \; 2 )  \\ 
 5,  &  j \equiv 1  \;  ({\rm mod} \; 2 ) \; {\rm and} \; i  = j   \; {\rm and} \; s \equiv  1 \;({\rm mod} \; 2 )  \\ 
 \end{array} \right.
\end{displaymath}
Analogously as in the proof of Theorem \ref{Pkn}, we can show that $f$ is a packing 5-coloring of $S_{C_{4t}}^n$.

Let $t\ge 2$, $k=4t+2$, $i,j \in [k]$ and $w \in [k]^{n-2}$. 
We define a coloring of  $wS_{C_{4t+2}}^{2}$ for every $w \in [k]^{n-2}$.
More precisely, we color $iC_{4t+2}$ of  $wS_{C_{4t+2}}^{2}$ such that we start at the vertex $ii$ and 
consecutively label its vertices with the  sequences presented in Table \ref{t4t+2}. 

\begin{table}[hbt]
\caption{A coloring of  $iC_{4t+2}$ in $S_{C_{4t+2}}^n$ } \label{t4t+2}
 \begin{tabular}{| c | c |c |c | } 
\multicolumn{4}{l}{$i < 4(t-1)$}  \\
 \hline
  $ i \equiv 0 \; ({\rm mod} \; 4 )$     & $ i \equiv 1 \; ({\rm mod} \; 4 )$ & $ i \equiv 2 \; ({\rm mod} \; 4)$ & $i \equiv  3  \;  ({\rm mod} \; 4 )$  \\
 \hline
 \\[-1em]
 $52142\overline{1312}1$ & $\overline{1213}2142131$ & $4\overline{1312}41231$ & $\overline{1312}412312$   \\ 
\hline 
\end{tabular}

\bigskip

 \begin{tabular}{|c | c | c |c |c | c |c | } 

\multicolumn{6}{l}{$i \ge 4(t-1)$}  \\
 \hline
$ \!\!i=4(t-1)\!\!$     &  $\!\! i=4(t-1)+1\!\!$   & $ \!\!i\!=\!\!4(t\!\!-\!\!1)\!\!+\!\!2\!\!$  & $\! i\!\!=\!\!4(t\!-\!1)\!+\!3\!\!$  & $\!\! i\!=\!4(t\!\!-\!\!1)\!\!+\!\!4\!\!$  & $ \!\!i\!\!=\!\!4(\!t\!-\!\!1)\!\!+\!\!5\!\!$  \\
 \hline
 \\[-1em]
 $\!\!51214\overline{1312}1\!\!$ & $\!\!14121\overline{3121}541231\!\!$ & $\!\!21\overline{3121}5131\!\!$ & $\!\!5\overline{1312}412315\!\!$  & $\!\!413\overline{1213}121\!\!$ & $\!\!\overline{1312}412312\!\!$  \\ 
\hline
\end{tabular}

\end{table}

In order to see that the coloring defined in Table  \ref{t4t+2} is a packing 5-coloring  of  $S_{C_{4t+2}}^n$, observe Fig. \ref{5C10} 
which represent the application of this  procedure to $S_{C_{10}}^2$. Since we can see that the obtained coloring is extendable, 
we showed that $ \chi_\rho(S_{C_{10}}^n)  \le  5$. Note that the procedure for $t \ge 3$ and $n=2$ repeats 
the coloring of the first four cycles (depicted by the dashed line in Fig. \ref{5C10}) in a way that the coloring remains extendable. 
Moreover, since the  coloring of each  $C_{4t+2}$ is extended either with the pattern $\overline{1312}$ or $\overline{1213}$, 
the conditions of extendable coloring 
are still fulfilled. Thus, we obtain an extendable  packing 5-coloring of $S_{C_{4t+2}}^2$ and the case is settled.

For  $k=4t+3$,  we define  a packing 5-coloring of $S_{C_k}^2$ as described in Table  \ref{t4t+3}. 
Since we can obtain by this procedure an extendable packing 5-coloring of $S_{C_{11}}^2$ depicted in Fig. \ref{5C11}, 
we can show by using the same arguments as above   that   $ \chi_\rho(S_{C_{4k+3}}^n)  \le  5$.

\begin{table}[hbt]
\caption{A coloring of  $iC_{4t+3}$ in $S_{C_{4t+3}}^n$ }  \label{t4t+3}
 \begin{tabular}{| c | c |c |c |c | } 
 \hline
  $\!\! i\! \equiv \!0 ({\rm mod}\,  4)\!\!$ & $\!\! i \!\equiv \! 1  ({\rm mod} \, 4 )\!\!$ & $\!\! i \!\equiv \! 2 ({\rm mod} \, 4),
   i\!\not=\!4t\!+\!2\!\!$ & $i \!\equiv \! 3   ({\rm mod} \, 4) \!\!$   &  $i=4t+2$ \\
 \hline
 \\[-1em]
 $\!\!1213124\overline{1213}1\!\!$ & $\!\!412134\overline{1213}1\!\!$ &  $1213124\overline{1213}1$ & $\!\!512134\overline{1213}1\!\!$  & $\!\!512134\overline{1213}2\!\!$ \\ 
\hline 
\end{tabular}
\end{table}

For  $k=4t+1$, $t \ge 3$,  we can construct a packing 5-coloring of $S_{C_k}^2$ from 
 an extendable  packing 5-coloring of $S_{C_{13}}^2$ depicted in Fig. \ref{5C13}. The construction 
for $t \ge 4$ repeats the first four cycles (depicted by the dashed line) of this graph as needed, while the coloring for each cycle is extended with the pattern
$\overline{1213}$ which can be found in every $C_{13}$ of  Figure \ref{5C13}. 
Note that the conditions of extendable coloring are still fulfilled with this procedure. 
Since we also found an extendable packing 5-coloring of $S_{C_9}^2$, 
it follows  that   $ \chi_\rho(S_{C_{4k+3}}^n)  \le  5$ for $t \ge 2$.
This assertion concludes the proof.
\end{proof}

\begin{figure}[hbt] 
\begin{center}

\includegraphics[width=120mm]{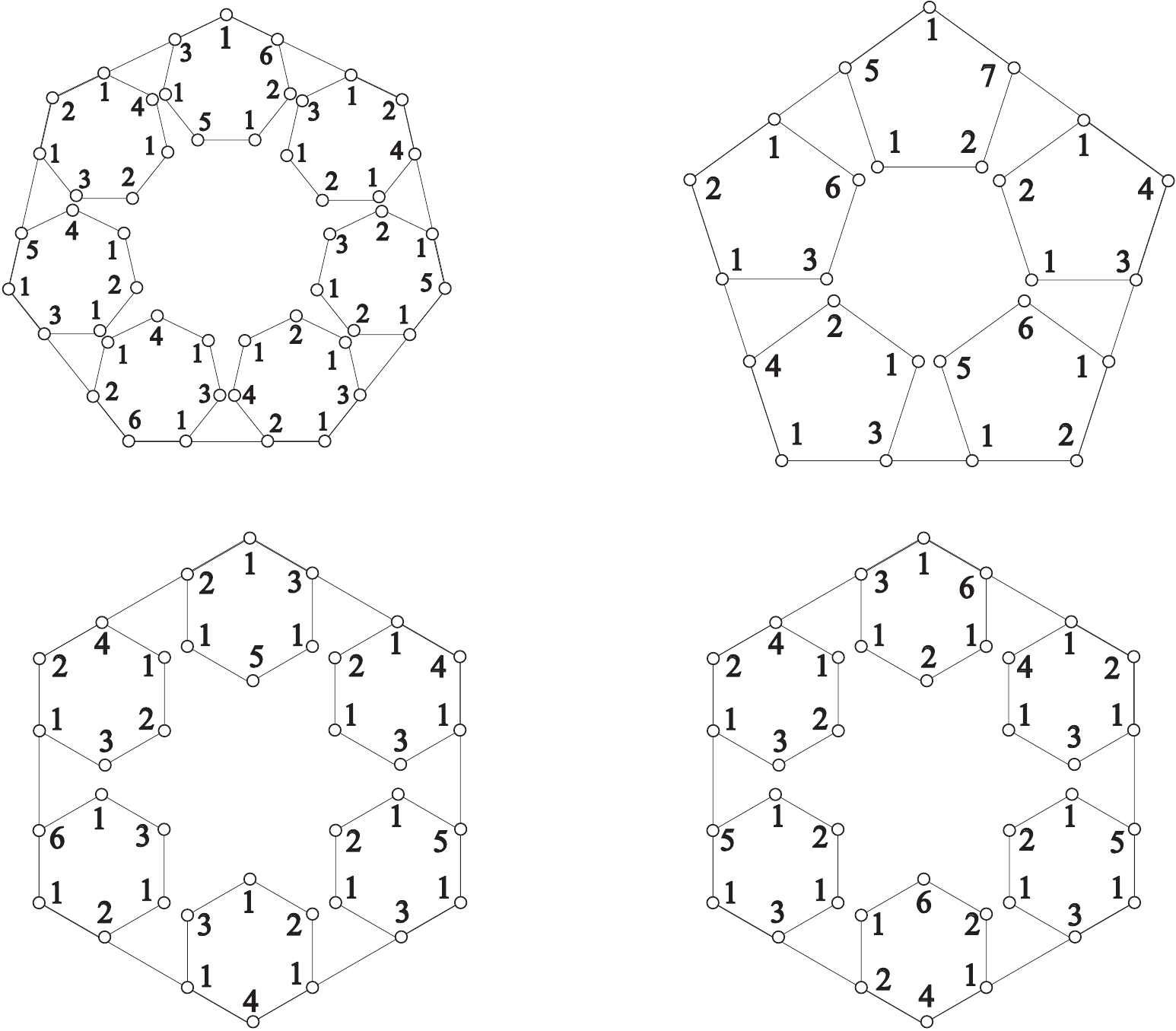} 

\caption{
 An extendable  packing 6-coloring of $S_{C_7}^2$,
 extendable  packing 7-coloring of $S_{C_5}^2$,  
 packing 5-coloring of $S_{C_6}^2$,   
 extendable  packing 6-coloring of $S_{C_6}^2$.
 } \label{cetvercek}
\end{center}
\end{figure}

We will need the following result provided by SAT solver.
\begin{prop1} \label{P6P8}
Let $S_{P_{k}}^2$ is the  generalized Sierpi\'nski  graph of  $P_{k}$ of dimension $2$.  
If $f$ (resp. $f'$) is an arbitrary packing 4-coloring of $S_{P_{6}}^2$  (resp. $S_{P_{8}}^2$), then  

(i) $f(22) \not \in \{2,3\}$ and

(ii) $f'(33)$ and $f'(44)$ are not both equal to 1.

\end{prop1}

\begin{cor1} \label{lower}
Let  $k$ be an odd integer. If $S_{C_{k}}^2$ is the  generalized Sierpi\'nski  graph of  $C_{k}$ of dimension $2$, then 
$\chi_\rho(S_{C_{k}}^2) \ge 5$.  
\end{cor1}

\begin{proof}
Let $k \ge 9$, $k$ is odd and $i \in [k]$. Suppose that  $f$ is a packing 4-coloring of $S_{C_{k}}^2$. 
Note that $S_{P_{6}}^2$  and $S_{P_{8}}^2$ are subgraphs of $S_{C_{k}}^2$. 
By Proposition \ref{P6P8}, for every external vertex $ii$ we have that $f(ii) \not \in \{2,3\}$. Moreover, 
two consecutive external vertices $ii$  and $(i+1)(i+1)$ (addition modulo $k$) cannot be both colored by 1.  
It follows that the external vertices of $S_{C_{k}}^2$
 are alternatively colored by 1 and 4. Since $k$ is odd, this is clearly impossible and we 
obtain a contradiction. It follows that  $\chi_\rho(S_{C_{k}}^2) \ge 5$ for every odd $k \ge 9$. 

Since we establish by SAT solver that a packing 4-coloring of $S_{C_{5}}^2$ and $S_{C_{7}}^2$ also do not exist, 
the proof is complete.
\end{proof}

Let $2\mathbb{N}$ and  $2\mathbb{N}+1$ denote the set of even and odd natural numbers, respectively.

\begin{thm1} \label{Ck}
Let $n \ge 2$ and $k \ge 4$. 
If $S_{C_{k}}^n$ is the  generalized Sierpi\'nski  graph of  $C_{k}$ of dimension $n$, then 
\begin{displaymath}
 \chi_\rho(S_{C_{k}}^n)  =  \left\{ \begin{array}{ll}
 4,  &    n = 2   \; {\rm and} \;  k \in 2\mathbb{N}\setminus \{6\}  \\
 5,  &   n \ge 3  \; {\rm and} \;   k \not\in \{5, 6, 7\} \; {\rm or} \;  n = 2  \; {\rm and} \;  k \in \{6\}   \cup 2\mathbb{N}+1\setminus \{5\}   \\
 6,  &   n \ge 3  \; {\rm and} \;   k \in \{ 6, 7\}  \; {\rm or} \; n \le 6  \; {\rm and} \;  k = 5    \\ 
  
 \end{array} \right.
\end{displaymath}
Moreover, if $n\ge 7$, then $6 \le  \chi_\rho(S_{C_{5}}^n)   \le 7$.

\end{thm1} 

\begin{proof}
First note that since $S_{P_{k}}^n$ is a subgraph of $S_{C_{k}}^n$, 
by Theorem \ref{Pkn} we have  $\chi_\rho(S_{C_{k}}^2)   \ge 4$ and  $\chi_\rho(S_{C_{k}}^n)   \ge 5$, $n\ge 3$. 
Moreover,  if $k$ is odd, we have $\chi_\rho(S_{C_{k}}^2) \ge 5$ by Corollary \ref{lower}.
Since we showed with  SAT solver that a packing 5-colorings of $S_{C_5}^2$,  
$S_{C_6}^3$ and $S_{C_7}^3$, as well as a packing 4-coloring of $S_{C_6}^2$ do not exist,
 all lower bounds are settled.

In order to prove the upper bounds, we found 
a packing 6-coloring of $S_{C_5}^6$,  
a packing 5-coloring of $S_{C_6}^2$,  
an extendable  packing 7-coloring of $S_{C_5}^2$,  
an extendable  packing 6-coloring of $S_{C_6}^2$ and
an extendable  packing 6-coloring of $S_{C_7}^2$ (the last four colorings are depicted in  Fig. \ref{cetvercek}).
Since the other needed upper bounds follows from Propositions  \ref{uC2} and \ref{uCk}, the proof is complete.

\end{proof}

\section{Sierpi\'nski  triangle graphs}

\begin{figure}[hbt] 
\begin{center}

\includegraphics[width=110mm]{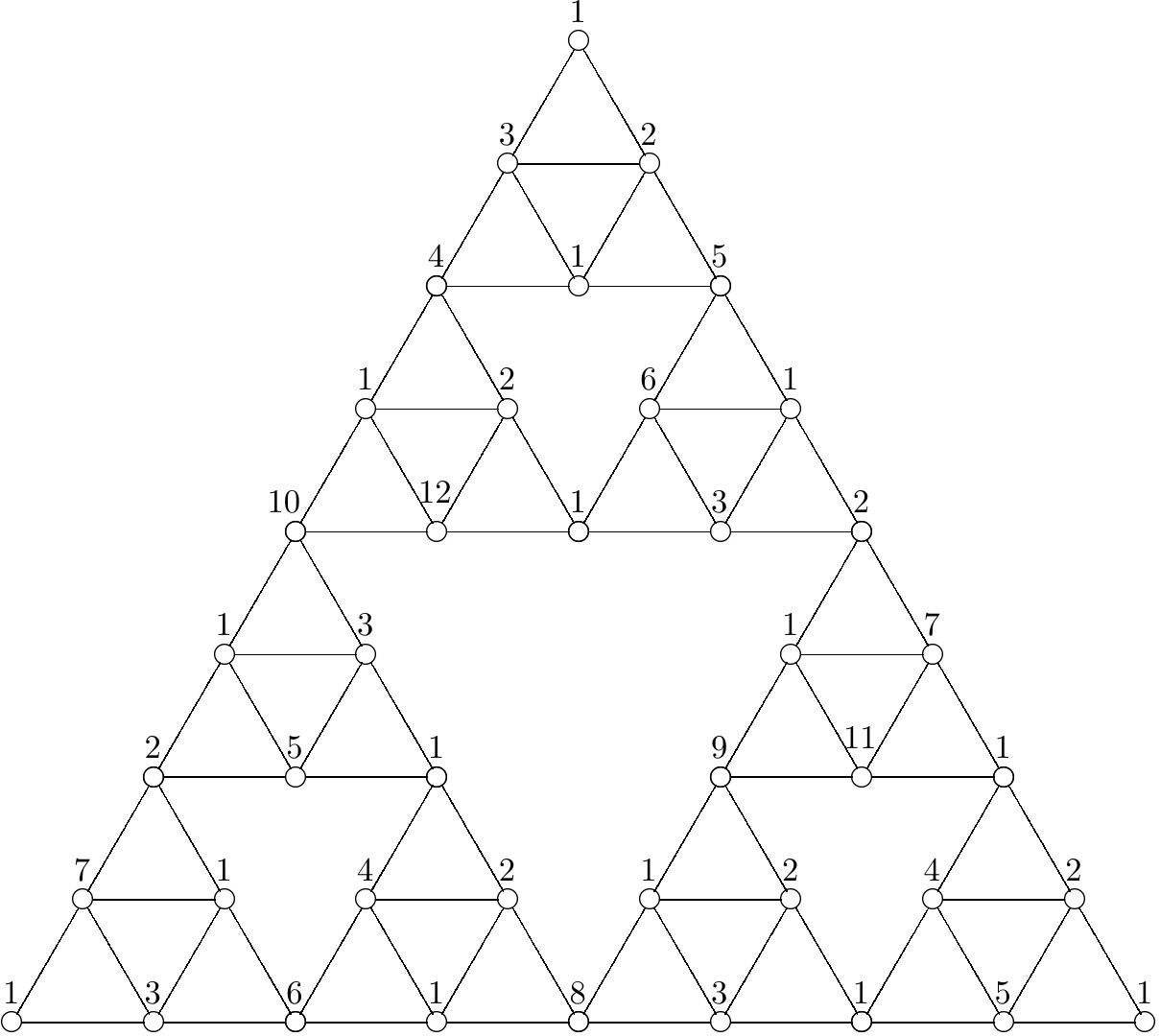} 

\caption{A packing 12-coloring of $ST_3^3$  } \label{12T3}
\end{center}
\end{figure}

The  {\em base-3 Sierpi\'{n}ski graphs $S^n$}  are defined such that we start with $S^0 = K_1$. 
For $n \ge 1$, the vertex set of $S^n$ is $[3]^n$ and the edge set is defined recursively as

$E(S^n) = \{\{is, i t\} : i \in [3], \{s, t\} \in E(S^{n-1})\} \cup
\{\{i j^{n-1}, ji^{n-1} \} | i, j \in [3], i\not= j \}$ .

As mentioned in the introduction, for $n\ge 1$, the base-3 Sierpi\'{n}ski graphs $S^n$ are generalized 
Sierpi\'{n}ski graphs  where $G=K_3$. 
Obviously, for $n\ge 1$,  $S^n$ can be constructed from three copies of  $S^{n-1}$. 

Let $n$ be a nonnegative integer. The class of the {\em Sierpi\'nski triangle graphs} $ST_3^n$ is obtained from 
$S^{n+1}$ by contracting all  non--clique edges.

\begin{figure}[hbt] 
\begin{center}

\includegraphics[width=130mm]{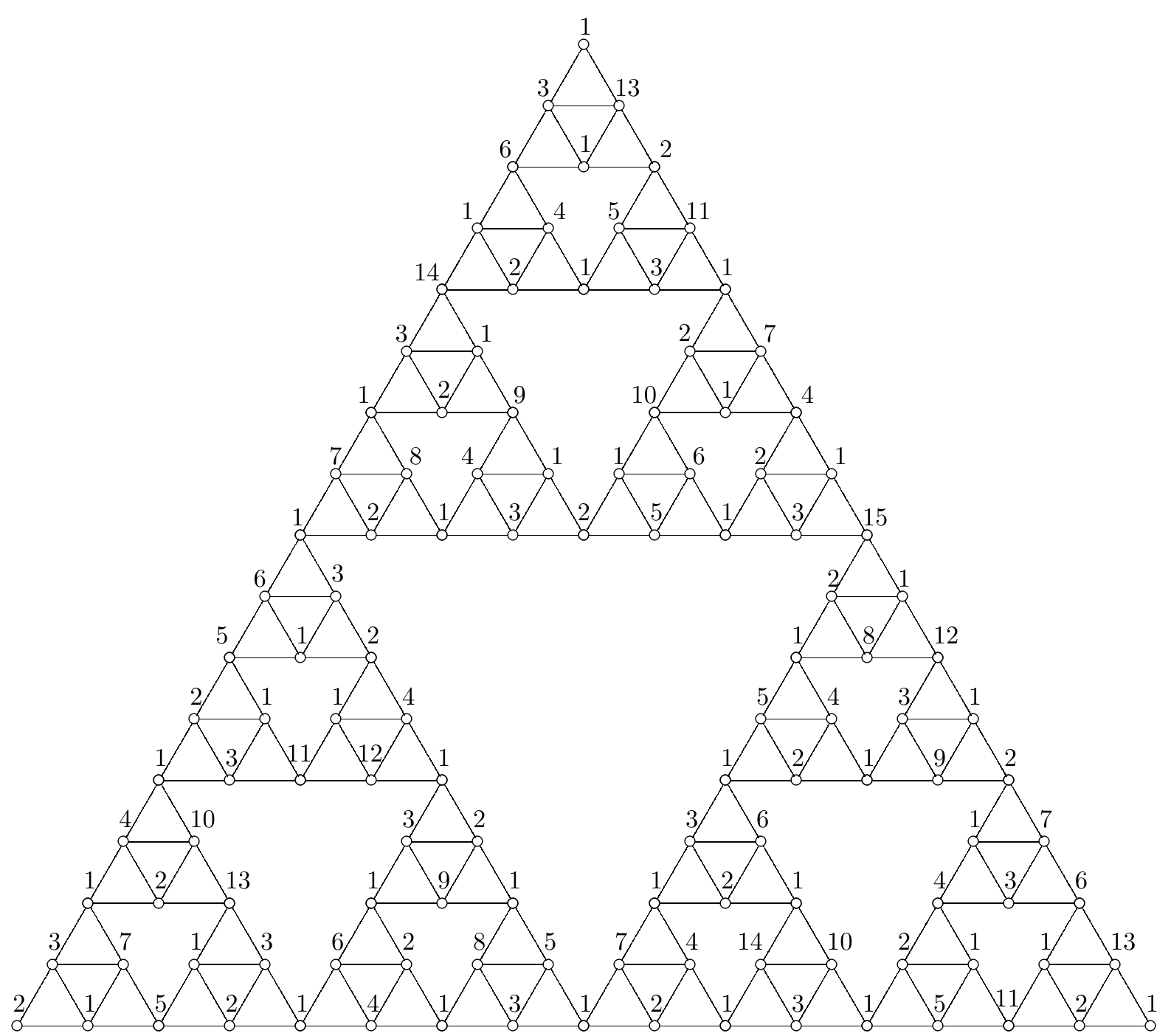} 

\caption{A packing 15-coloring of $ST_3^4$  } \label{15T4}
\end{center}
\end{figure}

There are various  other definitions of Sierpi\'nski triangle graphs, which are based on the fact that 
their drawings in the plane represent approximations of the Sierpi\'nski triangle fractal (see \cite{hinz}). 
More intuitively, 
Sierpi\'nski triangle graphs can be constructed by iteration. 
 We start with a complete graph on 3 vertices, i.e.\ $ST_3^0$ is the triangle $K_3$. 
Assume now  that $ST_3^{n}$ is already constructed. $ST_3^{n+1}$ is composed of three copies of $ST_3^{n}$, 
in a way that can be seen in Fig.  \ref{12T3}, where $ST_3^3$  is composed of three copies of $ST_3^2$.
Note that an extreme vertex of a copy of $ST_3^2$ is identified  with an extreme vertex of another copy 
(this procedure is done for exactly two extreme vertices of each copy).

\cite{BrFe} established the packing chromatic  number of Sierpi\'nski  triangle graphs  $ST_3^n$ for $n\le 2$ and 
showed that the packing chromatic number for this class of graphs  can be bounded above by 31. This bound is improved in the next theorem.

\begin{thm1}
Let $ST_3^n$ denote the Sierpi\'nski triangle graph of dimension $n$.
\begin{displaymath}
(i) \;  \chi_\rho(ST_3^n)  =  \left\{ \begin{array}{ll}
 3,  &  \;  n=0 \\ 
 4,  &  \;  n=1 \\ 
 8,  &  \;  n=2 \\ 
 12  &  \;  n=3. \\ 
 \end{array} \right.
\end{displaymath}

\begin{displaymath} (ii) \; 12 \le \chi_\rho(ST_3^4)  \le 15. \end{displaymath}
 
\begin{displaymath} (iii ) \; 12 \le \chi_\rho(ST_3^5)  \le 19. \end{displaymath}


Moreover, if $n \ge 6$, then  $ 12 \le \chi_\rho(ST_3^n)  \le 20$
\end{thm1}

\begin{proof}
For $n \le 2$, this is the result presented in \cite{BrFe}. 
Since we showed with  SAT solver that a packing 11-coloring of $ST_3^3$ cannot be obtained, we have 
$  \chi_\rho(ST_3^n) \ge 12 $ for every $n \ge 3$. 

The upper bounds were obtained 
by the following constructions: a packing 12-coloring of $ST_3^3$  presented in Fig. \ref{12T3}, 
a packing 15-coloring of $ST_3^4$  presented in Fig. \ref{15T4}, a packing 19-coloring of $ST_3^5$,   
while for $n \ge 6$ we obtained an extendable packing 20-coloring of  $ST_3^6$.
The last two constructions can be obtained from the authors or in the webpage 
\url{https://omr.fnm.um.si/wp-content/uploads/2017/06/SierpinskiP.pdf}.

\end{proof}


\nocite{*}
\bibliographystyle{abbrvnat}
\bibliography{articleSierpF}
\label{sec:biblio}

\end{document}